\documentclass{amsart}
\usepackage{ifthen,amsfonts,amsmath,amssymb,epic,eepic}
\usepackage{epsfig,color}

\newtheorem{theorem}{Theorem}[section]
\newtheorem{lemma}[theorem]{Lemma}
\newtheorem{proposition}[theorem]{Proposition}  
\newtheorem{corollary}[theorem]{Corollary}     

\theoremstyle{definition}

\theoremstyle{remark}
\newtheorem{remark}[theorem]{Remark}

\numberwithin{equation}{section}

\newcommand{\nn}{{\bf{n}}}
\newcommand{\Energy}{{\text{E}}}
\newcommand{\Ric}{{\text{Ric}}}

\newcommand{\Area}{{\text {Area}}}

\newcommand{\dist}{{\text {dist}}}
\newcommand\Tr{{\rm Tr}}

\newcommand{\K}{{\rm K}}
\newcommand{\cA}{{\mathcal{A}}}
\def\CC{{\bf  C}}
\def\RR{{\bf  R}}
\def\ZZ{{\bf  Z}}
\def\SS{{\bf  S}}

\newcommand{\eqr}[1]{(\ref{#1})}  

\begin{document}

\title[Extinction time for the Ricci flow and a question of Perelman]
{Estimates for the extinction time for the Ricci flow on certain
$3$--manifolds and a question of Perelman}

\author{Tobias H. Colding}%
\address{Courant Institute of Mathematical Sciences\\
251 Mercer Street\\ New York, NY 10012}
\author{William P. Minicozzi II}%
\address{Department of Mathematics\\
Johns Hopkins University\\
3400 N. Charles St.\\
Baltimore, MD 21218}
\thanks{The authors were partially supported by NSF Grants DMS
0104453 and DMS 0104187}

\email{colding@cims.nyu.edu and minicozz@math.jhu.edu}

\subjclass{Primary 53C44; Secondary 53C42, 57M50}



\keywords{Ricci flow, Finite extinction, $3$-manifolds, Min-max
surfaces}

\begin{abstract}
We show that the Ricci flow becomes extinct in finite time on any
Riemannian $3$--manifold without aspherical summands.
\end{abstract}

\maketitle

\section{introduction}

In this note we prove some bounds for the extinction time for the
Ricci flow on certain $3$--manifolds.  Our interest in this comes
from a question that Grisha Perelman asked the first author at a
dinner in New York City on April 25th of 2003.  His question was
``what happens to the Ricci flow on the $3$--sphere when one
starts with an arbitrary metric?  In particular, does the flow
become extinct in finite time?''  He then went on to say that one
of the difficulties in answering this is that he knew of no good
way of constructing minimal surfaces for such a metric in general.
However, there is a natural way of constructing such surfaces and
that comes from the min--max argument where the minimal of all
maximal slices of sweep--outs is a minimal surface; see, for
instance, \cite{CD}.  The idea is then to look at how the area of
this min--max surface changes under the flow.  Geometrically the
area measures a kind of width of the $3$--manifold and as we will
see for certain $3$--manifolds (those, like the $3$--sphere, whose
prime decomposition contains no aspherical factors) the area
becomes zero in finite time corresponding to that the solution
becomes extinct in finite time.
  Moreover, we will discuss a possible lower bound
for how fast the area becomes zero.  Very recently Perelman posted
a paper (see \cite{Pe1}) answering his original question about
finite extinction time.  However, even after the appearance of his
paper, then we still think that our slightly different approach
may be of interest.  In part because it is in some ways
geometrically more natural, in part because it also indicates that
lower bounds should hold, and in part because it avoids using the
curve shortening flow that he simultaneously with the Ricci flow
needed to invoke and thus our approach is in some respects
technically easier.

Let $M^3$ be a smooth closed orientable $3$--manifold and let
$g(t)$ be a one--parameter family of metrics on $M$ evolving by
the Ricci flow, so
\begin{equation}  \label{e:eqRic}
 \partial_t g=-2\,\Ric_{M_t}\, .
\end{equation}

Unless otherwise stated we will assume throughout that $M$ is
prime and non--aspherical (so $\pi_k (M) \ne \{ 0 \}$ for some
$k>1$). If $M$ is prime but not irreducible, then $M = \SS^2
\times \SS^1$ (proposition $1.4$ in \cite{Hr}) so $\pi_3 (M) =
\ZZ$. Otherwise, if $M$ is irreducible, then the sphere theorem
implies that   $\pi_2 (M) = 0$ (corollary $3.9$ in \cite{Hr}). In
the second case, the Hurewicz isomorphism theorem then implies
that $\pi_3 (M) \ne \{ 0 \}$ (since $M$ is non--aspherical).
Therefore, in either case, by suspension, as in lemma $3$ of
\cite{MiMo}, the space of maps from $\SS^2$ to $M$ is not simply
connected.

\begin{figure}[htbp]
\begin{center}
    \input{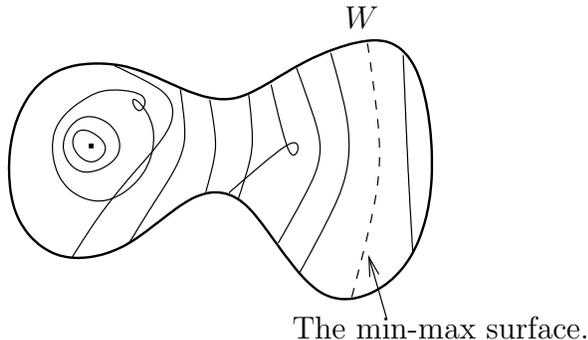}
    \caption{The sweep--out, the min--max surface, and the width W.}
    \label{f:1}
\end{center}
\end{figure}

Fix a continuous map $\beta: [0,1] \to C^0\cap L_1^2 (\SS^2 , M)$
where $\beta (0)$ and $\beta (1)$ are constant maps so that
$\beta$ is in  the nontrivial homotopy class $[\beta]$. We define
the width  $W=W(g,[\beta])$ by
\begin{equation}    \label{e:w3}
   W(g) = \min_{\gamma \in [ \beta]}
\, \max_{s \in [0,1]} \Energy (\gamma(s)) \,
   .
\end{equation}

One could equivalently define the width using the area rather than
the energy, but  the energy is somewhat easier to work with.  As
for the Plateau problem, this equivalence follows using the
uniformization theorem and the inequality $\Area (u) \leq \Energy
(u)$  (with equality when $u$ is a branched conformal map); cf.
lemma $4.12$ in \cite{CM}.\footnote{It may be of interest to
compare our notion of width, and the use of it, to a well--known
approach to the Poincar\'e conjecture.  This approach asks to show
that for any metric on a homotopy $3$--sphere a min--max type
argument produces an \underline{embedded} minimal $2$--sphere.
Note that in the definition of the width it plays no role whether
the minimal $2$--sphere is embedded or just immersed, and thus,
the analysis involved in this was settled a long time ago. This
well--known approach has been considered by many people, including
Freedman, Meeks, Pitts, Rubinstein, Schoen, Simon, Smith, and Yau;
see \cite{CD}.}

\vskip2mm The next theorem gives an upper bound for the derivative
of $W(g(t))$ under the Ricci flow which forces the solution $g(t)$
to become extinct in finite time (see paragraph 4.4 of \cite{Pe3}
for the precise definition of extinction time when surgery
occurs).

\begin{theorem}     \label{t:upper}
Let $M^3$ be a closed orientable prime non--aspherical
$3$--manifold equipped with a Riemannian metric $g=g(0)$. Under
the Ricci flow, the width $W(g(t))$ satisfies
\begin{equation}   \label{e:di1a}
\frac{d}{dt} \, W(g(t))  \leq -4 \pi + \frac{3}{4 (t+C)} \,
W(g(t))   \, ,
\end{equation}
in the sense of the limsup of forward difference quotients. Hence,
 $g(t)$  must become extinct in finite time.
\end{theorem}

The $4\pi$ in \eqr{e:di1a} comes from the Gauss--Bonnet theorem
and the $3/4$ comes from the bound on the minimum of the scalar
curvature that the evolution equation implies.  Both of these
constants matter whereas the constant $C$ depends on the initial
metric and the actual value is not important.

To see that \eqr{e:di1a} implies finite extinction time rewrite
\eqr{e:di1a} as
\begin{equation}
\frac{d}{dt} \left( W(g(t)) \, (t+C)^{-3/4} \right) \leq - 4\pi \,
(t+C)^{-3/4}
\end{equation} and
integrate to get
\begin{equation}  \label{e:lastaa}
 (T+C)^{-3/4} \, W(g(T)) \leq C^{-3/4} \, W(g(0))
- 16 \, \pi \, \left[ (T+C)^{1/4} - C^{1/4} \right]   \, .
\end{equation}
Since $W \geq 0$ by definition and the right hand side of
\eqr{e:lastaa} would become negative for $T$ sufficiently large we
get the claim.

\vskip2mm Arguing as in 1.5 of \cite{Pe1} (or alternatively using
Section \ref{s:4}), we get as a corollary of this theorem finite
extinction time for the Ricci flow on all $3$--manifolds without
aspherical summands.

\begin{corollary}  \label{c:upper}
Let $M^3$ be a closed orientable $3$--manifold whose prime
decomposition has only non--aspherical factors and is equipped
with a Riemannian metric $g=g(0)$. Under the Ricci flow with
surgery, $g(t)$ must become extinct in finite time.
\end{corollary}

Part of Perelman's interest in the question about finite time
extinction  comes from the following:  If one is interested in
geometrization of a homotopy three-sphere (or, more generally, a
three-manifold without aspherical summands) and knew that the
Ricci flow became extinct in finite time, then one would not need
to analyze what happens to the flow as time goes to infinity.
Thus, in particular, one would not need collapsing arguments.

\section{Upper bounds for the rate of change of area of
minimal $2$--spheres}

 Suppose that $\Sigma\subset M$ is a closed immersed
 surface (not necessarily minimal), then using \eqr{e:eqRic}
an easy calculation gives (cf. page 38--41 of \cite{Ha})
\begin{equation}        \label{e:diffAn}
\frac{d}{dt}_{t=0}\Area_{g(t)}( \Sigma ) =-\int_{\Sigma} [R -
\Ric_M (\nn,\nn)] \, .
\end{equation}
If $\Sigma$ is also minimal, then
\begin{align}        \label{e:diffA}
\frac{d}{dt}_{t=0}\Area_{g(t)}( \Sigma )
   &=-2\int_{\Sigma}\K_{\Sigma} -\int_{\Sigma}[|A|^2+\Ric_M
   (\nn,\nn)]  \\
  &=-\int_{\Sigma}K_{\Sigma}-\frac{1}{2}\int_{\Sigma}[ |A|^2 + R] \,
.   \notag
\end{align}
Here $\K_{\Sigma}$ is the (intrinsic) curvature of $\Sigma$, $\nn$
is a unit normal for $\Sigma$ (our $\Sigma$'s below will be
$\SS^2$'s and hence have a well--defined unit normal), $A$ is the
second fundamental form of $\Sigma$ so that $|A|^2$ is the sum of
the squares of the principal curvatures, $\Ric_M $ is the Ricci
curvature of $M$, and $R$ is the scalar curvature of $M$. (The
curvature is normalized so that on the unit $\SS^3$ the Ricci
curvature is $2$ and the scalar curvature is $6$.) To get
\eqr{e:diffA}, we used that by the Gauss equations and minimality
of $\Sigma$
\begin{equation}
\K_{\Sigma}=\K_M-\frac{1}{2}|A|^2\, ,
\end{equation}
where $\K_M$ is the sectional curvature of $M$ on the two--plane
tangent to $\Sigma$.

Our first lemma gives an upper bound for the rate of change of
area of minimal $2$--spheres.

\begin{lemma}     \label{l:upper}
If $\Sigma\subset M^3$ is a branched minimal immersion of the
$2$--sphere, then
\begin{equation}    \label{e:in24}
\frac{d}{dt}_{t=0}\Area_{g(t)}(\Sigma)  \leq -4 \pi -
\frac{\Area_{g(0)}(\Sigma)}{2} \, \min_{M} R(0) \, .
\end{equation}
\end{lemma}

\begin{proof}
Let $\{ p_i \}$ be the set of branch points of $\Sigma$ and $b_i >
0$ the order of branching at $p_i$.
 By \eqr{e:diffA}
\begin{equation}
 \frac{d}{dt}_{t=0} \Area_{g(t)}( \Sigma )
\leq -\int_{\Sigma}K_{\Sigma}-\frac{1}{2}\int_{\Sigma}R = -4\pi
-2\pi
 \sum b_i -\frac{1}{2}\int_{\Sigma}R  \, ,
\end{equation}
where the  equality used the Gauss--Bonnet theorem  with branch
points.  Note that branch points only help in the inequality
\eqr{e:in24}.
\end{proof}

 For an immersed minimal surface $\Sigma\subset M$ we set
\begin{equation}
    L_{\Sigma}\,\phi
    =\Delta_{\Sigma}\,\phi+|A|^2\,\phi+\Ric_M(\nn,\nn)\,\phi \, .
\end{equation}
 By \eqr{e:diffA} and the Gauss--Bonnet theorem
\begin{equation}   \label{e:gb}
\frac{d}{dt}_{t=0}\Area_{g(t)}(\Sigma)
=-2\int_{\Sigma}\K_{\Sigma}-\int_{\Sigma}1\cdot L_{\Sigma}\,1
=-4\pi\chi(\Sigma)-\int_{\Sigma}1\cdot L_{\Sigma}\,1\, .
\end{equation}
(Note that by the second variational formula (see, for instance,
section 1.7 of \cite{CM}), then
\begin{equation}
\frac{\partial^2}{\partial r^2}_{r=0}\Area (\Sigma_r)
=-\int_{\Sigma}\phi\,L_{\Sigma}\,\phi\, ,
\end{equation}
where $\Sigma_r=\{x+r\,\phi (x)\,\nn_{\Sigma} (x)\,|\,x\in
\Sigma\}$.)  Recall also that by definition the index of a minimal
surface $\Sigma$ is the number of negative eigenvalues (counted
with multiplicity) of $L_{\Sigma}$.  (A function $\eta$ is an
eigenfunction of $L_{\Sigma}$ with eigenvalue $\lambda$ if
$L_{\Sigma}\,\eta+\lambda\,\eta=0$.)   Thus in particular, since
$\Sigma$ is assumed to be closed, the index is always finite.

\section{Extinction in finite time}

 We begin by recalling a result on harmonic maps
which gives the existence of minimal spheres realizing  the width
$W(g)$.   The results of Sacks and Uhlenbeck give the harmonic
maps but potentially allow some loss of energy. This energy loss
was ruled out by Siu and Yau (using also arguments of Meeks and
Yau), see Chapter VIII in \cite{ScYa2}.  For our purposes, the
most convenient statement of this is given  in theorem $4.2.1$ of
\cite{Jo}.   (The bound for the index is not stated explicitly in
\cite{Jo} but follows immediately as in \cite{MiMo}.)

\begin{proposition}     \label{p:existence}
Given a metric $g$ on $M$ and a nontrivial $[\beta] \in \pi_1
(C^0\cap L^2_1 (\SS^2 , M))$, there exists  a sequence of
sweep--outs $\gamma^j: [0,1] \to C^0\cap L_1^2 (\SS^2 , M)$ with
$\gamma^j \in [\beta]$ so that
\begin{equation}
    W(g) = \lim_{j \to \infty} \, \max_{s \in [0,1]} \, \Energy
        (\gamma^j_s)   \, .
\end{equation}
Furthermore, there exist $s_j \in [0,1]$ and branched conformal
minimal immersions $u_0 , \dots , u_m : \SS^2 \to M$ with index at
most one so that, as $j \to \infty$, the maps $\gamma^j_{s_j}$
converge to $u_0$ weakly in $L_1^2$ and uniformly on compact
subsets of $\SS^2 \setminus \{ x_1 , \dots , x_k \}$, and
\begin{equation}
    W(g) =   \sum_{i=0}^m \Energy (u_i) = \lim_{j \to \infty} \Energy
        (\gamma^j_{s_j}) \,  .
\end{equation}
 Finally,
for each $i > 0$, there exists a point $x_{k_i}$ and a sequence of
conformal dilations $D_{i,j} : \SS^2 \to \SS^2$  about $x_{k_i}$
so that the maps $\gamma^j_{s_j} \circ D_{i,j}$ converge to $u_i$.
\end{proposition}

\begin{remark}
It is implicit in Proposition \ref{p:existence} that $W(g)
> 0$.  This can, for instance, be seen directly using \cite{Jo}.  Namely,
page 125 in \cite{Jo} shows that if $\max_s \Energy (\gamma^j_s)$
is sufficiently small (depending on $g$), then $\gamma^j$ is
homotopically trivial.
\end{remark}

 We will also need a  standard additional property for
 the min--max sequence of
  sweep--outs $\gamma^j$ of Proposition \ref{p:existence}
which can be achieved by modifying
  the sequence   as in section $4$ of \cite{CD}
(cf. proposition $4.1$ on page 85 in \cite{CD}).
  Loosely speaking this is the property that any subsequence
$\gamma_{s_k}^k$ with energy converging to $W(g)$ converges (after
possibly going to a further subsequence) to the union of branched
immersed minimal $2$--spheres each with index at most one.
Precisely this is that we can choose $\gamma^j$ so that:
 Given $\epsilon >
 0$, there exist $J$ and $\delta > 0$ (both depending on $g$ and $\gamma^j$)
so that if $j > J$ and
\begin{equation}    \label{e:eclose}
      \Energy (\gamma^j_s) > W(g) - \delta \, ,
\end{equation}
then there is a collection of branched minimal $2$--spheres $\{
\Sigma_i \}$ each of index at most one and with
\begin{equation}    \label{e:vclose}
    \dist  \, (\gamma^j_s ,
    \cup_i \Sigma_i ) < \epsilon \, .
\end{equation}
Here, the distance means  varifold distance (see, for instance,
section $4$ of \cite{CD}). Below we will use that, as an immediate
consequence of \eqr{e:vclose}, if $F$ is a quadratic form on $M$
and $\Gamma$ denotes $\gamma^j_s$, then
\begin{equation}     \label{e:hclose}
   \left| \int_{  \Gamma } [\Tr (F) - F(\nn_{\Gamma} ,
   \nn_{\Gamma})] - \sum_i \int_{\Sigma_i}
   [\Tr (F) - F(\nn_{\Sigma_i} ,
   \nn_{\Sigma_i})]  \right| <   C \, \epsilon \, \| F \|_{C^1}
\, \Area (\Gamma) \, .
\end{equation}

\vskip2mm
 In the proof of the result about finite extinction time we will also
need that the evolution equation for  $R= R(t)$, i.e. (see, for
instance, page 16 of \cite{Ha}),
\begin{equation}    \label{e:prescalar}
    \partial_t R = \Delta R + 2 |\Ric|^2   \geq \Delta R +
    \frac{2}{3} \, R^2 \, ,
\end{equation}
implies by a straightforward maximum principle argument that at
time $t > 0$
\begin{equation}    \label{e:scalar}
   R(t) \geq  \frac{1}{1/[ \min R(0)] -2t/3} =  -\frac{3}{2 (t+C)} \, .
\end{equation}
In the derivation of \eqr{e:scalar} we implicitly assumed that
$\min R(0)<0$.  If this was not the case, then \eqr{e:scalar}
trivially holds with $C=0$, since, by \eqr{e:prescalar}, $\min R
(t)$ is always non--decreasing.  This last remark is also used
when surgery occurs.  This is because by construction any surgery
region has large (positive) scalar curvature.

\vskip2mm
\begin{proof}
(of Theorem \ref{t:upper})   Fix a time $\tau$.  Below $\tilde{C}$
denotes a constant depending only on $\tau$ but will be allowed to
change from line to line.  Let $\gamma^j (\tau)$ be the sequence
of sweep--outs for the metric $g(\tau)$ given by
 Proposition \ref{p:existence}. We will use the sweep--out at
 time $\tau$ as a comparison to get an upper bound for the width
 at times $t > \tau$.   The key for this is the following claim
(the inequality
 in \eqr{e:acomp1} below):
Given $\epsilon > 0$, there exist $J$ and $\bar{h} > 0$ so that if
$j> J$ and $0 < h < \bar{h}$, then
\begin{align}        \label{e:acomp1}
    \Area_{g(\tau + h)}( \gamma^j_s (\tau ) )
    &- \max_{s} \, \Energy_{g(\tau)}( \gamma^j_s (\tau ))  \notag
    \\
&\leq
     [-4 \pi + \tilde{C} \, \epsilon +
 \frac{3}{4 (\tau+C)} \, \max_{s} \,
\Energy_{g(\tau)}( \gamma^j_s (\tau )) ] \, h + \tilde{C} \, h^2
\, .
\end{align}
 To see why
\eqr{e:acomp1} implies \eqr{e:di1a}, we use the definition of the
width to get
\begin{equation}        \label{e:defw}
    W(g (\tau + h) ) \leq \max_{s \in [0,1]}
    \Area_{g(\tau + h)}( \gamma^j_s (\tau ) ) \, ,
\end{equation}
and then  take the limit as $j\to \infty$ (so that $\max_{s} \,
\Energy_{g(\tau)}( \gamma^j_s (\tau ))\to W(g(\tau))$) in
\eqr{e:acomp1} to get
\begin{equation}        \label{e:defwq}
    \frac{W(g (\tau + h) ) - W(g (\tau ))}{h}
\leq    -4 \pi + \tilde{C} \, \epsilon +
 \frac{3}{4 (\tau+C)} \, W(g(\tau))     + \tilde{C} \, h
     \, .
\end{equation}
Taking $\epsilon \to 0$ in \eqr{e:defwq} gives \eqr{e:di1a}.

  It remains to prove
\eqr{e:acomp1}. First,  let $\delta > 0$ and $J$, depending on
$\epsilon$ (and on $\tau$), be given by
\eqr{e:eclose}--\eqr{e:hclose}. If $j > J$ and $\Energy_{g(\tau)}
(\gamma^j_s (\tau )  ) >  W(g) - \delta$, then let $\cup_i
\Sigma_{s,i}^j (\tau) $ be the collection of minimal spheres in
\eqr{e:hclose}. Combining \eqr{e:diffAn}, \eqr{e:hclose} with $F =
\Ric_M$, and Lemma \ref{l:upper} gives
\begin{align}        \label{e:diffAn2}
    \frac{d}{dt}_{t=\tau}\Area_{g(t)}( \gamma^j_s (\tau )) &
    \leq \frac{d}{dt}_{t=\tau}\Area_{g(t)}( \cup_i \Sigma_{s,i}^j (\tau) )
    +  \tilde{C} \, \epsilon \, \| \Ric_M \|_{C^1}
\, \Area_{g(t)}( \gamma^j_s (\tau ))\notag \\
        &\leq -4 \pi  - \frac{
\Energy_{g(\tau)}( \gamma^j_s (\tau ))}{2} \, \min_{M} R(\tau) +
\tilde{C} \, \epsilon \\
&\leq -4 \pi  +
 \frac{3}{4 (\tau+C)} \, \max_{s} \,
\Energy_{g(\tau)}( \gamma^j_s (\tau ))  + \tilde{C} \, \epsilon
\notag \, ,
\end{align}
where the last inequality used the lower bound \eqr{e:scalar} for
$R(\tau)$. Since the metrics $g(t)$ vary smoothly and every
sweep--out $\gamma^j$ has uniformly bounded energy, it is easy to
see that $\Energy_{g(\tau + h)} (\gamma^j_s (\tau )  )$ is a
smooth function of $h$ with a uniform $C^2$ bound independent of
both $j$ and $s$ near $h=0$ (cf. \eqr{e:diffAn}).  In particular,
\eqr{e:diffAn2} and Taylor expansion gives $\bar{h} > 0$
 (independent of $j$) so that \eqr{e:acomp1} holds
for  $s$ with $\Energy_{g(\tau)} (\gamma^j_s (\tau )  ) >  W(g) -
\delta$.  In
 the remaining case, we have  $\Energy (\gamma^j_s (\tau )) \leq
W(g) - \delta$ so the continuity of $g(t)$ implies that
\eqr{e:acomp1} automatically holds after possibly shrinking
$\bar{h}> 0$.

Finally, we claim that \eqr{e:di1a} implies finite extinction
time.   Namely, rewriting \eqr{e:di1a} as $\frac{d}{dt} \left(
W(g(t)) \, (t+C)^{-3/4} \right) \leq - 4\pi \, (t+C)^{-3/4}$ and
integrating gives
\begin{equation}  \label{e:last}
 (T+C)^{-3/4} \, W(g(T)) \leq C^{-3/4} \, W(g(0))
- 16 \, \pi \, \left[ (T+C)^{1/4} - C^{1/4} \right]   \, .
\end{equation}
Since $W \geq 0$ by definition and the right hand side of
\eqr{e:last} would become negative for $T$ sufficiently large, the
theorem follows.
\end{proof}

\section{Remarks on the reducible case} \label{s:4}

When $M$ is  reducible, then the factors in the prime
decomposition must split off in a uniformly bounded time. This
follows from a (easy) modification of the proof of Theorem
\ref{t:upper}. Namely, each (non--trivial) factor in the prime
decomposition gives rise to a $2$--sphere which does not bound a
$3$--ball in $M$ and, hence, to a stable minimal $2$--sphere in
this isotopy class by \cite{MeSiYa}. Applying the argument of the
proof of Theorem \ref{t:upper} to these  minimal $2$--spheres, we
see that the minimal area in this isotopy class must go to zero in
finite time as claimed.

\appendix

\section{Lower bounds for the rate of change of area of
minimal $2$--spheres}

 The next lower bound is an adaptation of Hersch's
theorem; cf. \cite{ChYa}. Recall that Hersch's theorem (see, for
instance, \cite{ScYa1}) states the sharp scale invariant
inequality that for any metric on the $2$--sphere $\lambda_1$
times the area is bounded uniformly from above by the
corresponding quantity on a round $2$--sphere.

\begin{lemma}  \label{l:l1}
If $\Sigma\subset M^3$ is an immersed minimal $2$--sphere with
index at most one, then
\begin{equation}  \label{e:el1}
8\,\pi\geq \int_{\Sigma}[|A|^2+\Ric_M(\nn,\nn)] =
\int_{\Sigma}1\cdot L_{\Sigma}\,1\, .
\end{equation}
Hence, by \eqr{e:gb}
\begin{equation}
\frac{d}{dt}_{t=0}\Area_{g(t)}(\Sigma) \geq - 16\, \pi \, .
\end{equation}
\end{lemma}

\begin{proof}
If $\Sigma$ is stable, i.e., if the index is zero, then for all
$\phi$
\begin{equation}
-\int_{\Sigma}\phi\,L\,\phi\geq 0\, ,
\end{equation}
or equivalently
\begin{equation}   \label{e:vardesc}
\int_{\Sigma}|\nabla \phi|^2\geq
\int_{\Sigma}[|A|^2+\Ric_M(\nn,\nn)]\,\phi^2\, ,
\end{equation}
and thus by letting $\phi\equiv 1$ in \eqr{e:vardesc} we see that
\eqr{e:el1} holds.

If the index is one, then we let $\eta$ be an eigenfunction for
$L_{\Sigma}$ with negative eigenvalue $\lambda<0$.  That is,
\begin{equation} \label{e:defpsi}
L_{\Sigma} \,\eta+\lambda\,\eta=0\, .
\end{equation}
By a standard argument, then an eigenfunction corresponding to the
first eigenvalue of a Schr\"odinger operator (Laplacian plus
potential) does not change sign and thus we may assume that $\eta$
is everywhere positive.  In particular, $\int_{\Sigma}\eta>0$.
Since $\Sigma$ has index one then it follows that \eqr{e:vardesc}
holds for all $\phi$ with
\begin{equation}
0=\int_{\Sigma}\eta\, \phi\, .
\end{equation}

By the uniformization theorem, there exists a conformal
diffeomorphism $\Phi:\Sigma\to \SS^2\subset \RR^3$.   For
$i=1,2,3$ set $\phi_i=x_i\circ \Phi$.  For $x\in \SS^2$ let
$\pi_x:\SS^2\setminus \{x\}\to \CC$ be the stereographic
projection and let $\psi_{x,t}(y)=\pi_x^{-1}(t(\pi_x(y)))$, then
for each $t,\,x$ this can be extended to a conformal map on
$\SS^2$. Define $\Psi:\SS^2\times [0,1)\to G$, where $G$ is the
group of conformal transformations of $\SS^2$, by $\Psi
(x,t)=\psi_{x,1/(1-t)}$.  Since $\Psi (x,0)=\text{id}_{\SS^2}$ for
each $x\in \SS^2$, $\Psi$ can be thought of as a continuous map on
$B_1(0)=\SS^2\times [0,1)/(x,0)\equiv (y,0)$.  Set
\begin{equation}
\cA (\Psi (x,t))=\frac{1}{\int_{\Sigma} \eta} \left(
\int_{\Sigma}\eta\, x_i\circ\Psi (x,t)\circ\Phi \right)
_{i=1,2,3}\, ,
\end{equation}
where $\eta$ is as in \eqr{e:defpsi}.  It follows that
\begin{equation}
\cA:B_1(0)\to B_1(0) \hbox{ and } \lim_{(y,t)\to
(x,1)}\cA(\Psi(y,t))= x \, .
\end{equation}
In particular, it follows that $\cA$ extends to $\partial B_1(0)$
as the identity map.  We can therefore, by elementary topology
(after possibly replacing $\Phi$ by $\psi\circ\Phi$), assume that
for each $i$
\begin{equation}  \label{e:e2}
\int_{\Sigma}\eta\,\phi_i=0\, ;
\end{equation}
that is each $\phi_i$ is orthogonal to $\eta$. It follows from
\eqr{e:vardesc} that for each $i$
\begin{equation}  \label{e:e3}
\int_{\Sigma}|\nabla \phi_i|^2\geq
\int_{\Sigma}[|A|^2+\Ric_M(\nn,\nn)]\,\phi^2_i\, .
\end{equation}
Summing over $i$ and using that $\Phi(\Sigma)\subset \SS^2$ so
$\sum_{i=1}^3\phi_i^2=1$ we get
\begin{equation}  \label{e:e4}
\sum_{i=1}^3\int_{\Sigma}|\nabla \phi_i|^2 \geq
\int_{\Sigma}[|A|^2+\Ric_M(\nn,\nn)]\, .
\end{equation}
Now obviously, since $\Phi$ is conformal (so that it preserves
energy) and since each $x_i$ is an eigenfunction for the Laplacian
on $\SS^2\subset \RR^3$ with eigenvalue $\lambda_1(\SS^2)=2$, we
get
\begin{equation}  \label{e:e5}
\int_{\Sigma}|\nabla \phi_i|^2=\int_{\SS^2}|\nabla x_i|^2
=\lambda_1 (\SS^2)\int_{\SS^2}x_i^2\, .
\end{equation}
Combining \eqr{e:e4} with \eqr{e:e5} we get
\begin{align}  \label{e:e6}
2\,\Area (\SS^2)&=\sum_{i=1}^3\lambda_1(\SS^2)\int_{\SS^2}x_i^2
=\sum_{i=1}^3\int_{\SS^2}|\nabla x_i|^2\\
&=\sum_{i=1}^3\int_{\Sigma}|\nabla \phi_i|^2 \geq
\int_{\Sigma}[|A|^2+\Ric_M(\nn,\nn)]\, .
\end{align}
\end{proof}

\vskip2mm {\it Acknowledgement: We are grateful to John Lott for
helpful comments.}

\bibliographystyle{amsplain}

\end{document}